\documentclass{amsart}

\usepackage{amsmath,amsfonts,amssymb,amsthm}

\def\refpp#1{(\ref {#1})}
\newcommand{\mspan}{\mathrm{span}}

\def\labtt#1{\label {#1}  }

\def\a{\alpha}
\def\b{\beta}
\def\g{\gamma}

\def\o{\omega}

\def\CC{{\mathbb C}}

\def\RR{{\mathbb R}}

\def\ZZ{{\mathbb Z}}

\def\la{\langle}
\def\ra{\rangle}
\def\<{\langle}
\def\>{\rangle}

\def\half{{1 \over 2}}
\def\fourth{{1 \over 4}}

\def\sixteenth{{1 \over 16}}

\def\dual#1{#1^*}        

\def\allone{\hbox{\bf 1}}

\newcommand{\dg}[1]{\mathcal{D}(#1)}
\begin{document}

\newtheorem{thm}{Theorem}[section]
\newtheorem{prop}[thm]{Proposition}
\newtheorem{lem}[thm]{Lemma}
\newtheorem{rem}[thm]{Remark}
\newtheorem{coro}[thm]{Corollary}
\newtheorem{conj}[thm]{Conjecture}
\newtheorem{de}[thm]{Definition}
\newtheorem{hyp}[thm]{Hypothesis}

\newtheorem{nota}[thm]{Notation}
\newtheorem{ex}[thm]{Example}
\newtheorem{proc}[thm]{Procedure}

\title{Rootless pairs of $EE_8$-lattices}
\author{Robert L.~Griess, Jr.}
\address[R.\,L.~Griess, Jr.]{
Department of Mathematics, University of Michigan,
Ann Arbor, MI 48109-1043,
USA}
\email{rlg@umich.edu}

\author{Ching Hung Lam}
\address[C.\,H. Lam]{
Department of Mathematics and
National Center for Theoretical Sciences,
National Cheng Kung University,
Tainan, Taiwan 701}

\email{chlam@mail.ncku.edu.tw}

\begin{abstract}
We describe a classification of
pairs $M, N$ of lattices isometric to $EE_8  := \sqrt 2  E_8$ such that the lattice $M + N$ is integral and rootless and such that the dihedral group associated to them has order at most 12. It turns out that most of these pairs may be embedded in the Leech lattice.  Complete proofs will appear in another article.  This theory of integral lattices has connections to vertex operator algebra theory and moonshine.
\end{abstract}

\maketitle


\section{Introduction}
Associated to a sublattice $S$ of an integral lattice $L$ is an orthogonal
transformation of order 2 on Euclidean space (defined to be $-1$ on $S$ and 1 on the annihilator of $S$). This involution leaves  $L$ invariant if $S$ satisfies the \textit{RSSD condition} ($2L \le S + ann(S )$) \cite{bwy, cbwy}.   
An example of such an $S$ is $\sqrt 2 U$, where $U$ is an integral unimodular lattice, such as $EE_8:=\sqrt 2 E_8$.  
A pair of RSSD sublattices gives a dihedral group in the isometry group of $L$.  

In this article, we describe our classification of  pairs of $EE_8$-lattices
which span an integral and rootless lattice and whose associated involutions
(isometries of order 2) generate a dihedral group of order at most
12. 
There are exactly $11$ such pairs.  Most can be realized inside the Leech lattice, which we denote by $\Lambda$.  
We refer to an upcoming article  \cite{gl} for proofs and details.

This work may be considered purely  as a study of positive definite integral lattices.
Our real motivation, however, is the evolving theory of vertex operator algebras (VOA) and their automorphism groups.\footnote{AMS classifications: 20C10 (integral representations of finite groups); 11H55 (quadratic forms); 11F22 (relations to Lie algebras and finite groups).}

\section{Motivation from VOA theory}
The primary connection between the Monster and vertex operator
algebras was established in \cite{flm}.
In \cite{Mi}, Miyamoto gave a theory for defining involutions on a VOA.
Suppose that $(V, Y, \o, \allone)$ is a VOA. A degree $2$ element $e\in V$ is called a \textit{rational conformal vector of central charge $1/2$} if the
components of $Y(e,z)=\sum_{n\in \ZZ} e_n z^{-n-1}$
(The $e_i$ are in $End(V)$)
together with $Id_V$ span a Virasoro Lie algebra with central charge $\half$ and the subVOA of $V$ generated by $w$ is a simple VOA, so is isomorphic to simple Virasoro VOA $L(\half, 0)$.

The simple Virasoro VOA $L(\half, 0)$ has only three irreducible modules,
$L(\half, 0)$, $L(\half, \half )$, $L(\half, \frac 1{16})$, and the fusion rules give two formulas for an involution on $V$.

Since $L(\frac{1}{2},0)$ is rational, i.e., all $L(\frac{1}{2},0)$-modules are completely reducible, we can decompose $V:=V_0 \oplus V_{\half } \oplus V_{\sixteenth}$ into the direct sum of $V_c$, where
$V_c$ denotes the sum of all $L(\frac{1}{2},0)$-submodules of $V$ isomorphic to $L(\half, c)$.

\medskip

(1)  Define a linear map  $\tau_e:V \to V$ to be $1$ on $V_0 \oplus V_\half$ and $-1$ on $V_\sixteenth$. Miyamoto\,\cite{Mi} showed that $\tau_e$ is always an automorphism of $V$. This automorphism is called \textit{a Miyamoto involution} if $\tau_e$ has order $2$.

\medskip

(2) If $\tau_e=Id_V$, i.e., $V_{\sixteenth}=0$, then the linear map $\sigma_e$ defined by $1$ on $V_0$ and $-1$ on $V_\half$ is also an automorphism. It is called \textit{a Miyamoto involution of $\sigma$-type}.

\medskip

These involutions have been studied a lot.
Miyamoto proved \cite{Mi1}
the important result
that there is a bijection between the conjugacy class of $2A$
involutions in the Monster simple group and conformal vectors of
central charge $\half$ in the moonshine vertex operator algebra
$V^\natural$.

\subsection{Lattice type vertex operator algebras}
Recall that a lattice VOA $V_L$ associated with an even positive definite lattice $L$ is given by $V_L=S(\hat{\mathfrak{h}}^-) \otimes \CC [L]$, as a vector space, where
 $\mathfrak{h}:=\CC \otimes L$, $\hat{\mathfrak{h}}^-= \oplus_{n=1}^\infty t^{-n} \otimes \mathfrak{h}$,  $S(\hat{\mathfrak{h}}^-)$ is the symmetric algebra of
$\hat{\mathfrak{h}}^-$ and $\CC[L]=\mathrm{span}_\CC\{ e^\a| \a\in L\}$ is the group algebra of $L$.

Let $V_L$ be a lattice VOA.  It is common to study an involution on
$V$ which is a ``lift'' of the $-1$ map on $L$.
The standard lift of $-1$, called $\theta$, acts as $-1$ on $\hat{\mathfrak{h}}^-$ and takes
 $e^\a$ to $e^{-\a}$, for all $\a$.
The fixed point subVOA of $\theta$ is denoted $V_L^+$.

For VOAs of lattice type $V_L$ or $V_L^+$, there are two known explicit general formulas for rational conformal vectors of central charge $1/2$.

\medskip

\noindent {\it First formula: }

When $\a$ is a norm 4 vector of $L$,
\begin{equation}\labtt{aa1cv}
e=\sixteenth (t^{-1}\otimes \a)^2 \pm \fourth (e^\a +e^{-\a})
\end{equation}
is a rational conformal vector of central charge $1/2$.

\medskip

\noindent {\it Second formula: }

When $E$ is a sublattice of $L$, $E\cong EE_8$,
\begin{equation}\labtt{ee8cv}
e=\frac 18 q +\frac 1{32}\sum_{\a \in E/\{\pm 1\} } \varphi (\a) (e^\a + e^{-\a})
\end{equation}
is also  a rational conformal vector of central charge $1/2$,
(here $q$ is a sum $\sum_i  (t^{-1}\otimes u_i)^2$, where $u_i$ is an orthonormal basis of  $\CC\otimes E$;
$\varphi$ is a homomorphism $E\rightarrow \{\pm 1\}$).

\medskip

This second formula explains our interest in $EE_8$.

\subsection{Why rootless?}

If $L$ is rootless, $V_L^+$ has zero degree 1 term and so (by general VOA theory) the degree 2 part is a commutative algebra (generally nonassociative).  In this case, the automorphism group of $V_L^+$ is finite (not generally true for VOAs). Moreover, it is conjectured \cite{lsy} that the two kinds of conformal vectors in \refpp{aa1cv} and \refpp{ee8cv} will exhaust all the rational conformal vectors of
central charge $1/2$ in $V_L^+$ if $L$ is rootless. This conjecture was proved when $L$
is a $\sqrt{2}$ times a root lattice or the Leech lattice
\cite{lsy,ls} but it is still open if $L$ is a general rootless
lattice.   The results of this paper could help settle this
conjecture.

\subsection{Why dihedral groups of order at most 12? }
We may consider pairs of Miyamoto involutions and what they generate.  For a general VOA, this is an arbitrary dihedral group.
In the case of $V_L^+$ when $L$ is rootless,  the possibilities are limited to dihedral of order at most 12, by a recent theorem of Shinya Sakuma \cite{Sa}.  In fact, his
hypotheses are more general: the VOA need not be lattice type, but it should be defined over $\RR$ with a positive definite form, have no negative terms, the degree 0 component should be 1-dimensional, and degree 1 component should be 0-dimensional.  His proof is quite technical.

Even without the Sakuma theorem, one knows that
the restriction to
dihedral groups of order at most $12$ is relevant to a study of sporadic groups.  Consider the 2A class of involutions in the Monster \cite{gms} and classes within the Baby Monster, $Fi_{24}$, etc., discussed in \cite{gn}.

\section{The classification of integral, rootless $EE_8$-pairs}

In this section, we shall discuss our classification of $EE_8$-pairs. The following is our main theorem.

\begin{thm} Suppose that in Euclidean space, the integral lattice $L$ is the sum $M+N$ of sublattices, where $M\cong N \cong EE_8$, and $L$ is rootless. Let $t_M$ be the involution defined by $M$ and let $t_N$ be the  involution defined by $N$.  Assume that the dihedral group $\la t_M, t_N \ra$ has order at most 12.
Then the possibilities
for $L$ are
listed in Table \ref{NREE8SUM} and all these possibilities exist.  The lattices in Table \ref{NREE8SUM} are unique up to isometry of pairs $M, N$.  Except for $DIH_4(15)$, all of them embed as sublattices of the Leech lattice.
\end{thm}

\def\dih#1#2{DIH_{#1}(#2)}

\begin{table}[bht]
\caption{ \bf  NREE8SUMs: integral rootless lattices which are sums of
$EE_8$s}
\begin{center}
\begin{tabular}{|c|c|l|c|c|}
\hline
 Name & $\<t_M,t_N\>$ &Isometry type of $L$\,(contains)& $\dg L$ & In\,Leech?\cr
 \hline \hline
 $\dih{4}{12}$ & $Dih_4$ & $\geq  DD_4^{\perp 3}$ & $1^4 2^6 4^2$ & Yes\cr
 \hline
 $\dih{4}{14}$ & $Dih_4$ & $\geq AA_1^{\perp 2} \perp DD_6^{\perp 2}$ & $1^4 2^8 4^2$  &Yes\cr
 \hline
 $\dih{4}{15}$ & $Dih_4$ & $\geq AA_1\perp EE_7^{\perp 2}$ & $1^22^{14}  $  & No\cr
 \hline $\dih{4}{16}$ & $Dih_4$ & $\cong EE_8\perp EE_8$ & $2^{16}$ &   Yes\cr
 \hline $\dih{6}{14}$ & $Dih_6$ & $\geq AA_2\perp A_2\otimes E_6 $ & $1^7 3^3 6^2$ &Yes\cr
 \hline
 $\dih{6}{16}$ & $Dih_6$ & $\cong A_2\otimes E_8 $ & $1^8 3^8 $ &Yes\cr
 \hline
 $\dih{8}{15}$ & $Dih_8$ & $ \geq AA_1^{\perp 7}\perp EE_8$     & $1^{10}4^5 $&   Yes\cr
 \hline
 $\dih{8}{16, DD_4}$ & $Dih_8$ & $\geq DD_4^{\perp 2} \perp  EE_8$  & $1^8 2^4 4^4$  &Yes\cr
 \hline
 $\dih{8}{16, 0}$ & $Dih_8$ & $\cong BW_{16} $ & $1^8 2^8$ &Yes\cr
 \hline $\dih{10}{16}$ & $Dih_{10}$ & $\geq A_4\otimes A_4 $ & $1^{12}5^4$  &Yes\cr
 \hline $\dih{12}{16}$ & $Dih_{12}$ & $\geq AA_2\perp AA_2 $ & &\cr
 & & $\ \perp A_2\otimes E_6 $ & $1^{12}6^4 $  & Yes\cr
 \hline \hline
\end{tabular}
\medskip

$X^{\perp n}$ denotes the orthogonal sum of $n$ copies of the lattice $X$.
\end{center}
\label{NREE8SUM}
\end{table}%

\begin{nota}
$A_1, \dots, E_8$ denote the root lattices of the corresponding type;

$AA_1$, $\dots$, $EE_8$ denote the lattices isometric to $\sqrt{2}$ times the lattices $A_1, \dots, E_8$;

$\dg L = \dual L / L$ is the discriminant group of $L$;

$\dih {r}s$ means a pair of $EE_8$ lattices $(M, N)$ such that $t_M$, $t_N$ generates a dihedral group of order $r$ and $rank(M+N)=s$;

$DIH_8(16, P)$ means $DIH_8(16)$ and that $P= ann_M(N)\cong ann_N(M)$.
\end{nota}

\subsection{Some highlights of the proof}\labtt{hlight}
Next, we shall discuss the main steps for the classification.
We go through cases $|t_Mt_N| = 2,3,4,5,6$.  Our respective analyses are called {\it $DIH_4$-theory, $DIH_6$-theory,
$DIH_{8}$-theory, $DIH_{10}$-theory,  $DIH_{12}$-theory.  }

In $\la t_M, t_N \ra$, let $g:=t_Mt_N$.   Then $\ZZ [\la g \ra]$ acts on $L$ and it acts on $J:=ann_L(Fix_L(g))$, where $Fix_L(g)$ denotes the set of all fixed points of $g$ in $L$.  The action is that of as a ring of integers in a number field when $|g|$ is prime.  

The main idea is to determine possibilities for $Fix_L(g)$, $J$, $ann_M(N)$, $ann_N(M)$ and related sublattices.
Exhaustive case by case analysis gives a list of candidates.  In all cases, the candidates are proved unique, given certain things we deduce about their sublattices.

\medskip

First, we observe that $\la t_M, t_N \ra$ acts faithfully on $L$ and leaves invariant  $Fix_L(g)=M\cap N$.

\medskip

\noindent \textit{Case:} $|g|=2,3,5$.

In these cases,  we  determine all sublattices of $E_8$ which are direct summands and whose discriminant group  is an elementary abelian $p$-group for $p=2, 3, 5$.  Exhaustive case by case analysis gives a list of candidates. It turns out  $M\cap N$ is isometric to $\sqrt 2$ times one of these lattices.  In fact,   $M\cap N\cong 0, AA_1, AA_1\perp AA_1$, or $DD_4$ if $|g|=2$; $M\cap N\cong 0$ or $AA_2$ if $|g|=3$; and $M\cap N=0$ if $|g|=5$.

\medskip
Given $M\cap N$, we then analyse $J$ and its sublattices.

\medskip

When $|g|=2$, $\la t_M, t_N \ra$ is a four-group. Then  we have $M\cap J=ann_M(N)$, $N\cap J=ann_N(M)$ and $M\cap N \perp ann_M(N)\perp ann_N(M)$ is an index $2$ sublattice of $L$.  In this case, the isometry type of $L$ is uniquely determined by $M\cap N$.

\medskip

When $|g|=3$, we consider the $\ZZ [\la g \ra]$-submodule $K$ generated by $M\cap J$.  Then $K$ is a sublattice of $J$ and $K$ is isomorphic as a lattice to $A_2\otimes \frac 1{\sqrt 2}(M\cap J)$.  The possibilities for $M\cap J$ in this case are $EE6$ or $EE8$. Again, the isometry type of $L$  is uniquely determined by $M\cap N$.

\medskip

When $|g|=5$, $M\cap N=0$. We show that for any norm $4$ vector $\a\in L$, the $\ZZ [\la g \ra]$-submodule generated by $\a$ is isomorphic as a lattice to $AA_4$. In fact, we show that $L=M+N$ contains a sublattice $U$ isometric to the orthogonal sum of $4$ copies of $AA_4$ such that $M\cap U\cong N\cap U \cong AA_1^8$. The uniqueness of $L$ is then shown by explicit glueings.

\medskip

\noindent \textit{Case:} $|g|=4,6$.

In each of these cases, let $h:=g^2$. Then $(M, Mh)$ and $(N, Nh)$ are $EE_8$ pairs whose associated dihedral group has order $4$ or $6$. We then use the results for $Dih_4$ and $Dih_6$ to deduce the structures of $L$. It turns out that there is only one possible case for $|g|=6$ but $3$ different cases for $|g|=4$.

A proof that the candidates are really rootless is made easier by a magic tool.  Most candidates can be embedded in the Leech lattice by direct constructions (and use of a uniqueness result).    Since Leech is rootless, so is our candidate $L$.  See the last section for examples of embeddings.

\section{ Glauberman-Norton subgroups}

In \cite{gn},  Glauberman-Norton computed the centralizers in the Monster of a pair of 2A-involutions.
We hope to study these subgroups by promoting pairs $M, N$  of lattices from our list to involutions on the Moonshine VOA.
An interesting point is that the Glauberman-Norton subgroups  have ``Levi factors'' which look like index $2$ subgroups of a direct product of a pair of Weyl groups modulo centers.  Our theory may give a context for this by viewing such a subgroup as acting on the Moonshine VOA.   At this time, we can see actions of direct products of Weyl groups on large subVOAs of the Moonshine VOA,
in the following way.
The Moonshine VOA is a twisted lattice type VOA based on the Leech lattice.  The Leech lattice contains sublattices of types $AA_1$ or $EE_8$ which will give a concrete realization of each pair of $2A$-involutions in the Glauberman-Norton theory (by Formula \refpp{aa1cv} or Formula \refpp{ee8cv}).  Since there are examples of tensor products of lattices contained in our   $M+N$ when $|g|=3$ or $5$ \refpp{hlight}, we have actions of a finite group with quotient a direct product of Weyl groups
on a subVOA of lattice type.

\section{Examples of some pairs $M, N$ in the Leech lattice}

In this section,
we shall give a few examples of $EE_8$-pairs $M, N$ in the Leech lattice.
Because of space limitations, these are just for the $DIH_6$-theory.
First, we recall some notation for describing the Leech lattice and its isometries \cite{cs, Gr12}.

\begin{nota}
Arrange  the set $\Omega=\{1, \dots, 24\}$ into a $4\times 6$ array such that the six columns forms a sextet.
For each codeword in the Golay code $\mathcal{G}$, $0$ and $1$ are marked by a
blanked and non-blanked space, respectively, at the corresponding
positions in the array. For example, $(1^8 0^{16})$ is denoted by
the array
\[
\begin{array}{|cc|cc|cc|}
\hline *& *&\ &\ &\ &\ \\
* & *&\ &\ &\ &\ \\
* & *&\ &\ &\ &\ \\
* & *&\ &\ &\ &\ \\ \hline
\end{array}
\]

By the same notation, every vector in the Leech lattice $\Lambda$ can be written as the form
\[
X=\frac{1}{\sqrt{8}}\left[ X_1 X_2 X_3 X_4 X_5 X_6\right],\quad \text{
juxtaposition of column vectors}.
\]
For example,
\[
\frac{1}{\sqrt{8}}\,
\begin{array}{|rr|rr|rr|}
 \hline 2  & 2& 0  & 0 & 0 & 0 \\
 2 &  2 &  0&  0& 0& 0\\
 2 &  2& 0& 0&  0 & 0 \\
 2 &  2& 0&  0& 0 & 0\\ \hline
\end{array}.
\]
\end{nota}

\begin{de}\labtt{defofxi}

Let
\[
A= \frac{1}2 \left [
\begin{matrix}
-1 & 1& 1& 1\\
1 & -1& 1& 1\\
1 & 1& -1& 1\\
1 & 1& 1& -1
\end{matrix}\right].
\]

Define a linear map $\xi:\Lambda \to \Lambda$  by  $ X\xi = AXD$, where
\[
X=\frac{1}{\sqrt{8}}\left[ X_1 X_2 X_3 X_4 X_5 X_6\right]
\]
is a vector in the Leech lattice $\Lambda$ and $D$ is the diagonal matrix
\[
\begin{pmatrix}
-1&0&0&0&0&0\\
0& 1&0&0&0&0\\
0& 0& 1&0&0&0\\
0&0&0 &1&0&0\\
0&0&0&0&1&0\\
0&0&0&0&0& 1\\
\end{pmatrix}.
\]
Recall that $\xi$ defines an isometry of $\Lambda$ (cf. \cite[p. 288]{cs} and \cite[p. 97]{Gr12}).

\end{de}

\begin{nota}
 Let $\mathcal{O}$ be an octad. Denote by $E(\mathcal{O})$ the sublattice of the Leech lattice $\Lambda$ which is supported on $\mathcal{O}$. Then $E(\mathcal{O})\cong EE_8$.
\end{nota}

Here are examples for the $DIH_6$ cases.

\subsubsection{$ |g|=3$.}\labtt{g=3} In this case, $M\cap N =0$ or
$AA_2$.

\medskip
\noindent $\bf \dih{6}{16}$

Let
\[
\mathcal{O}_1=
\begin{array}{|cc|cc|cc|}
\hline  * & * &  &  &  & \\
* & * &\ &\ &\ &\ \\
* & * &\ &\ &\ &\ \\
* & * &\ &\ &\ &\ \\ \hline
\end{array}
\]
be an octad and $M:=E(\mathcal{O}_1)\cong EE_8$. We choose a basis
$\{\b_1,\dots, \b_8\}$ of $M$, where
\begin{align*}
\b_1= \frac{1}{\sqrt{8}}\,
\begin{array}{|cc|cc|cc|}
\hline 4 & -4& 0  & 0 & 0 & 0 \\
0 &  0&  0&  0& 0& 0\\
 0 & 0&  0& 0&  0 & 0 \\
 0 &  0& 0&  0& 0 & 0\\ \hline
\end{array}, &&
 \b_2= \frac{1}{\sqrt{8}}\,
 \begin{array}{|rr|rr|rr|}
\hline 0  & 4& 0  & 0 & 0 & 0 \\
      -4 &  0&  0&  0& 0& 0\\
       0 & 0&  0& 0&  0 & 0 \\
       0 &  0& 0&  0& 0 & 0\\ \hline
\end{array}, \\
\b_3= \frac{1}{\sqrt{8}}\,
 \begin{array}{|rr|rr|rr|}
\hline 0  & 0& 0  & 0 & 0 & 0 \\
       4 &  -4&  0&  0& 0& 0\\
       0 & 0&  0& 0&  0 & 0 \\
       0 &  0& 0&  0& 0 & 0\\ \hline
\end{array}, &&
\b_4=\frac{1}{\sqrt{8}}\,
\begin{array}{|rr|rr|rr|}
\hline 0  & 0& 0  & 0 & 0 & 0 \\
       0&  4 &  0&  0& 0& 0\\
       -4 & 0&  0& 0&  0 & 0 \\
       0 &  0& 0&  0& 0 & 0\\ \hline
\end{array},
\\
\b_5=\frac{1}{\sqrt{8}}\,
\begin{array}{|rr|rr|rr|}
\hline 0  & 0& 0  & 0 & 0 & 0 \\
       0&  0 &  0&  0& 0& 0\\
       4 & -4&  0& 0&  0 & 0 \\
       0 &  0& 0&  0& 0 & 0\\ \hline
\end{array},  &&
\b_6=\frac{1}{\sqrt{8}}\,
\begin{array}{|rr|rr|rr|}
\hline 0  & 0& 0  & 0 & 0 & 0 \\
       0 & 0&  0&  0& 0& 0\\
       0 & 4&  0& 0&  0 & 0 \\
       -4&  0& 0&  0& 0 & 0\\ \hline
\end{array}, \\
\b_7=\frac{1}{\sqrt{8}}\,
\begin{array}{|rr|rr|rr|}
\hline -4  & -4& 0  & 0 & 0 & 0 \\
       0&  0 &  0&  0& 0& 0\\
       0 & 0&  0& 0&  0 & 0 \\
       0 &  0& 0&  0& 0 & 0\\ \hline
\end{array}, &&
\b_8=\frac{1}{\sqrt{8}}\,
\begin{array}{|rr|rr|rr|}
\hline\ 2  &\ 2& 0  & 0 & 0 & 0 \\
       2 &  2 &  0&  0& 0& 0\\
       2 &  2& 0& 0&  0 & 0 \\
       2 &  2& 0&  0& 0 & 0\\ \hline
\end{array}.
\end{align*}

Let $N$ be the lattice generated by the vectors
\begin{align*}
&\a_1=\frac{1}{\sqrt{8}}\,
\begin{array}{|rr|rr|rr|}
\hline\ \ 2 & -2&\ \ 2  & -2 &\ \ 2 &\ \, 0 \\
\ 0 & \ 0&\ 0&\ 0&\ 0&  2\\
\ 0 &  \ 0&\ 0&\ 0&\ 0 & 2 \\
\ 0 & \ 0&\ 0&\ 0&\ 0 &\ \, 2\\ \hline
\end{array}, &&
 \a_2=\frac{1}{\sqrt{8}}\,
\begin{array}{|rr|rr|rr|}
\hline 0 &\ 2&\ 0&\ 2&  0 &\ 2 \\
      -2 &\ 0& -2&\ 0&\ 0 &\ -2 \\
     \ 0 &\ 0&\ 0&\ 0&\ 2 &\ -2\\
     \ 0 &\ 0&\ 0&\ 0&\ 0 &\ 0 \\ \hline
\end{array},\\
&\a_3=\frac{1}{\sqrt{8}}\,
\begin{array}{|rr|rr|rr|}
\hline \ 0 &\  0&\ 0  &\ 0 &\ 0 & -2 \\
     \ \ 2 & -2 & \ \,2 & -2 & 2  &\ 0\\
       \ 0 &\ 0 &\ 0  &\ 0 & -2  &  0\\
       \ 0 &\ 0 &\ 0  &\ 0 & 2  &\ 0 \\ \hline
\end{array}, &
&\a_4=\frac{1}{\sqrt{8}}\,
\begin{array}{|rr|rr|rr|}
\hline  \ 0 &\ 0  &   \ 0&    0&\ 0 &\ 2 \\
        \ 0 &\ \ 2&   \ 0&\ \ 2&  -2&\ \ 2\\
         -2 &\   0&    -2&\   0&\ 0 & \ 0\\
        \ 0 &  \ 0&   \ 0& \  0&\ 0 & -2 \\ \hline
\end{array},\\
&\a_5=\frac{1}{\sqrt{8}}\,
\begin{array}{|rr|rr|rr|}
\hline   \ 0 &\ 0&  \ 0&\ 0& \ 0 &-2 \\
         \ 0 &\ 0&  \ 0&\ 0&\ 2& 0\\
       \ \ 2 & -2&\ \ 2& -2&\ 2 &0\\
         \ 0 &\ 0&\ 0&  0& -2 & 0 \\ \hline
\end{array}, &
&\a_6=\frac{1}{\sqrt{8}}\,
\begin{array}{|rr|rr|rr|}
\hline \ 0 &  \ 0&\ 0 &  \ 0& 0   &\ 2 \\
       \ 0 &  \ 0& \ 0&  \ 0& 0   & -2\\
       \ 0 &\ \ 2&   0&\ \ 2&\ -2   &\ 2\\
        -2 &  \ 0&  -2&    0&\ 0 &\ \ 0 \\ \hline
\end{array},\\
&\a_7=\frac{1}{\sqrt{8}}\,
\begin{array}{|rr|rr|rr|}
\hline     -2 &  -2&  -2& -2& \ \ 0 & -2 \\
           \ 0 &\ 0& \ 0&\ 0& -2& \ 0\\
           \ 0 &\ 0& \ 0&\ 0& -2 &\ 0\\
           \ 0 &\ 0&   0&\ 0& -2 &\ 0 \\ \hline
\end{array}, &
&\a_8=\frac{1}{\sqrt{8}}\,
\begin{array}{|rr|rr|rr|}
\hline \ \ 1 &\ \ \, 1&\ \ 1  &\ \ \,1& \ -3 & 1 \\
       \ \ 1 &\ \ 1&\ \ 1  &\ \ 1& \  \ 1 & 1 \\
       \ \ 1 &\ \ 1&\ \ 1  &\ \ 1& \  \ 1 & 1 \\
       \ \ 1 &\ \ 1&\ \ 1  &\ \ 1& \  \ 1 & 1 \\ \hline
\end{array}.
\end{align*}
By checking the inner products, it is easy to shows that $N\cong EE_8$. Note that $\a_1, \dots, \a_7 $ are supported on octads and thus $N\leq \Lambda$.

In this case, $M\cap N=0$ and $\{\b_1, \b_2, \dots, \b_8, \a_1,\dots, \a_8\}$ is a basis of $L=M+N$. The Gram matrix of $L=M+N$ is given by

\[
\left[
\begin{array}{rrrr rrrr| rrrr rrrr}
4&-2&0&0&0&0&0&0&2&-1&0&0&0&0&0&0\\
-2&4&-2&0&0&0&-2&0&-1&2&-1&0&0&0&-1&0\\
0&-2&4&-2&0&0&0&0&0&-1&2&-1&0&0&0&0\\
0&0&-2&4&-2&0&0&0&0&0&-1&2&-1&0&0&0\\
0&0&0&-2&4&-2&0&0&0&0&0&-1&2&-1&0&0\\
0&0&0&0&-2&4&0&0&0&0&0&0&-1&2&0&0\\
0&-2&0&0&0&0&4&-2&0&-1&0&0&0&0&2&-1\\
0&0&0&0&0&0&-2&4&0&0&0&0&0&0&-1&2\\ \hline
2&-1&0&0&0&0&0&0&4&-2&0&0&0&0&0&0\\
-1&2&-1&0&0&0&-1&0&-2&4&-2&0&0&0&-2&0\\
0&-1&2&-1&0&0&0&0&0&-2&4&-2&0&0&0&0\\
0&0&-1&2&-1&0&0&0&0&0&-2&4&-2&0&0&0\\
0&0&0&-1&2&-1&0&0&0&0&0&-2&4&-2&0&0\\
0&0&0&0&-1&2&0&0&0&0&0&0&-2&4&0&0\\
0&-1&0&0&0&0&2&-1&0&-2&0&0&0&0&4&-2\\
0&0&0&0&0&0&-1&2&0&0&0&0&0&0&-2&4
\end{array}
\right]
\]

\medskip

By looking at the Gram matrix, it is clear that $L=M+N\cong A_2\otimes E_8$. The Smith invariant
sequence is $ 1111 1111 3333 3333$.

\medskip

\noindent $\bf \dih{6}{12}$

Let
\[
\mathcal{O}_2=\begin{array}{|cc|cc|cc|}
\hline \ & *&* &* &* &* \\
* &\ &\ &\ &\ &\ \\
* &\ &\ &\ &\ &\ \\
* &\ &\ &\ &\ &\ \\ \hline
\end{array}
\]
and  $M:=E(\mathcal{O}_2)$. Set $N:=M\xi $, where  $\xi$ is the
isometry defined in  \refpp{defofxi}. Consider
\begin{align*}
\g_1= \frac{1}{\sqrt{8}}\,
\begin{array}{|cc|cc|cc|}
\hline 0 & 0& 0  & 0 & 0 & 0 \\
 0 &  0&  0&  0& 0& 0\\
 -4 & 0&  0& 0&  0 & 0 \\
 4 &  0& 0&  0& 0 & 0\\ \hline
\end{array}, &&
 \g_2= \frac{1}{\sqrt{8}}\,
 \begin{array}{|rr|rr|rr|}
\hline 0  & 0& 0  & 0 & 0 & 0 \\
      -4 &  0&  0&  0& 0& 0\\
       4 & 0&  0& 0&  0 & 0 \\
       0 &  0& 0&  0& 0 & 0\\ \hline
\end{array}, \\
\g_3= \frac{1}{\sqrt{8}}\,
 \begin{array}{|rr|rr|rr|}
\hline 0  & -4& 0  & 0 & 0 & 0 \\
       4 &  0&  0&  0& 0& 0\\
       0 & 0&  0& 0&  0 & 0 \\
       0 &  0& 0&  0& 0 & 0\\ \hline
\end{array}, &&
\g_4=\frac{1}{\sqrt{8}}\,
\begin{array}{|rr|rr|rr|}
\hline 0  & 4& -4  & 0 & 0 & 0 \\
       0 &  0 &  0&  0& 0& 0\\
       0 & 0&  0& 0&  0 & 0 \\
       0 &  0& 0&  0& 0 & 0\\ \hline
\end{array},\\
\g_5=\frac{1}{\sqrt{8}}\,
\begin{array}{|rr|rr|rr|}
\hline 0  & 0& 4  & -4 & 0 & 0 \\
       0 & 0&  0&  0& 0& 0\\
       0 & 0&  0& 0&  0 & 0 \\
       0 &  0& 0&  0& 0 & 0\\ \hline
\end{array}, &&
\g_6=\frac{1}{\sqrt{8}}\,
\begin{array}{|rr|rr|rr|}
\hline 0  & 0& 0  & 4 & -4 & 0 \\
       0&  0 &  0&  0& 0& 0\\
       0 & 0&  0& 0&  0 & 0 \\
       0 &  0& 0&  0& 0 & 0\\ \hline
\end{array}, \\
\g_7=\frac{1}{\sqrt{8}}\,
\begin{array}{|rr|rr|rr|}
\hline 0  & 0& 0  & 0 & 0 & 0 \\
       0&  0 &  0&  0& 0& 0\\
       -4 & 0&  0& 0&  0 & 0 \\
       -4 &  0& 0&  0& 0 & 0\\ \hline
\end{array}, &&
\g_8=\frac{1}{\sqrt{8}}\,
\begin{array}{|rr|rr|rr|}
\hline\ 0  &\ 2& 2  & 2 & 2 & 2 \\
       2 &  0 &  0&  0& 0& 0\\
       2 &  0& 0& 0&  0 & 0 \\
       2 &  0& 0&  0& 0 & 0\\ \hline
\end{array}.
\end{align*}
Then $\{\g_1, \dots,\g_8\}$ is a basis of $M$ and $\{\g_1\xi ,
\dots,\g_8\xi \}$ is a basis of $N$.

\medskip

By the definition of $\xi$, it is easy to show that $\g_1\xi=-\g_1,
\g_2\xi=-\g_2$. Moreover, for any $\a\in M=E(\mathcal{O}_2)$, $\a
\xi$ is supported on $\mathcal{O}_2$ if and only if $\a\in
\mspan_\ZZ\{ \g_1, \g_2\}$.  Hence, $F=M\cap N=\mspan_\ZZ\{ \g_1,
\g_2\} \cong AA_2$. Then $ann_M(F)\cong ann_N(F)\cong EE_6$ and
$L=M+N$ is of rank $14$.

Note that
$\{\g_1, \g_2, \g_3, \dots, \g_8\}$ is a basis and $\{\g_1, \g_2, \g_3\xi, \dots, \g_8\xi \} $ is a basis of $N$. Therefore,
$$\{\g_1, \g_2\} \cup \{ \g_3, \dots, \g_8\}\cup \{\g_3\xi, \dots, \g_8\xi \} $$
is a basis of $L$ and the Gram matrix of $L$ is given by

\[
\left[
\begin{array}{cccccccccccccc}
  4&-2&0&0&0&0&0&0&0&0&0&0&0&0\\
  -2&4&-2&0&0&0&-2&0&-2&0&0&0&2&0\\
  0&-2&4&-2&0&0&0&0&0&1&0&0&-2&0\\
  0&0&-2&4&-2&0&0&0&1&-2&1&0&0&0\\
  0&0&0&-2&4&-2&0&0&0&1&-2&1&0&0\\
  0&0&0&0&-2&4&0&0&0&0&1&-2&0&0\\
  0&-2&0&0&0&0&4&-2&2&0&0&0&0&1\\
  0&0&0&0&0&0&-2&4&0&0&0&0&-1&-2\\
  0&-2&0&1&0&0&2&0&4&-2&0&0&0&0\\
  0&0&1&-2&1&0&0&0&-2&4&-2&0&0&0\\
  0&0&0&1&-2&1&0&0&0&-2&4&-2&0&0\\
  0&0&0&0&1&-2&0&0&0&0&-2&4&0&0\\
  0&2&-2&0&0&0&0&-1&0&0&0&0&4&2\\
  0&0&0&0&0&0&1&-2&0&0&0&0&2&4
  \end{array}
  \right]
\]
whose Smith invariant sequence is $1111\ 11111 \ 333 66$.

In this case,  $ann_M(F)+ann_N(F) \cong A_2\otimes E_6$ and thus $L$ contains a sublattice isometric to $AA_2 \perp (A_2\otimes E_6)$.

{}

\end{document}